\documentclass{article}
\usepackage{spconf,amsmath,graphicx}
\usepackage{xcolor}
\usepackage{amssymb}
\usepackage{graphicx}
\usepackage{subfig}
\usepackage{mwe} % to get dummy images
\usepackage{algorithm}
\usepackage{algorithmic}
\usepackage{dsfont}
\usepackage{hyperref}
\usepackage{pgfplots,tikz}
\usetikzlibrary{arrows,shapes,shadows.blur,spy,fadings,shadings}
\usepackage{makecell}
\usepackage{diagbox}
\usepackage{changes}
\usepackage{cancel}

\pgfplotsset{compat=newest}
\pgfplotsset{plot coordinates/math parser=false}
\newlength\figureheight
\newlength\figurewidth

\def\lspace{$\ell_0$\;}	% l0_
\def\CL{COL0RME}
\def\Phie{\Phi_{{\scriptscriptstyle \mathtt{CEL0}}}} % CEl0 penalty in ND
\DeclareMathOperator{\EX}{\mathbb{E}}% expected value

\definecolor{myred}{RGB}{150,20,20}

\DeclareMathOperator*{\argmin}{arg\,min}

\title{\CL: COvariance-based \lspace super-Resolution Microscopy with intensity Estimation}
 \name{\begin{tabular}{c}Vasiliki Stergiopoulou$^{1}$, José Henrique de Morais Goulart$^{2}$, Sébastien Schaub$^{3}$\\
Luca Calatroni$^{1}$, and Laure Blanc-Féraud$^{1}$\end{tabular}}
\address{$^{1}$ Université Côte d’Azur, CNRS, INRIA, I3S, France \\
$^{2}$ IRIT, Université de Toulouse, CNRS, Toulouse INP, Toulouse, France \\
$^{3}$ Sorbonne Université, CNRS, LBDV, France \\}

\begin{document}
%\ninept
%
\maketitle
\begin{abstract}
Super-resolution light microscopy overcomes the physical barriers due to light diffraction, allowing the observation of otherwise indistinguishable subcellular entities. %However, under challenging experimental conditions, state-of-the art super-resolution methods, sometimes, do not achieve adequate spatio-temporal resolution.
However, the specific acquisition conditions required by state-of-the-art super-resolution methods to achieve adequate spatio-temporal resolution are often very demanding. Exploiting molecules fluctuations allows good spatio-temporal resolution live-cell imaging by means of common microscopes and conventional fluorescent dyes.  
In this work, we present the method \CL\ for COvariance-based \lspace super-Resolution Microscopy with intensity Estimation. It codifies the assumption of sparse distribution of the fluorescent molecules as well as the temporal and spatial independence between emitters via a non-convex optimization problem formulated in the covariance domain. In order to deal with real data, the proposed approach also estimates background and noise statistics. It also includes a final estimation step where intensity information is retrieved, which is valuable for biological interpretation and future applications to super-resolution imaging.

\end{abstract}
\vspace{0.07cm}
\section{Introduction}
\label{sec:intro}

Fluorescence microscopy is an imaging technique that allows the investigation of living cells and their organelles. However, due to light diffraction, structures smaller than half the emission wavelength (typically, around 250 nm) cannot be resolved. As many entities of biological interest have a size smaller than such a barrier, it is therefore fundamental to consider approaches which are able to provide a super-resolved version of the acquired data. Techniques such as Single Molecule Localization Microscopy (SMLM) \cite{smlm} and STimulated Emission Depletion (STED) \cite{Hell:94} can achieve nanometric spatial resolution at the cost of low temporal resolution and need of special (often costly) equipment. On the other side, Structured Illumination Microscopy (SIM) \cite{sim} has fast acquisition rates but relatively low spatial acquisition.  

% On the other hand, Structured Illumination Microscopy (SIM) \cite{sim} allows for high speed acquisitions using standard dies and cheaper special equipment, but reaches relatively low spatial resolutions compared to the two methods above.
% should i mention sparsity-based methods (e.g. FALCON) ? 
% Laure : I don't think so. Falcon is one the SMLM methods,and if we cite Falcon we have to cite also the other approaches. We just quote SMLM, it is enough. 
%Vasiina: ok :)
To overcome these limitations, different super-resolution approaches taking advantage of the independent stochastic temporal fluctuations of \emph{standard} fluorescent dyes were considered during the last decade. These methods take as an input a temporal stack of images acquired by common microscopes and produce a super-resolved image on a finer grid. In Super-resolution Optical Fluctuation Imaging (SOFI) \cite{sofi}, for instance, fluctuations are exploited by mapping the acquired image data into the covariance domain. While SOFI can drastically improve the temporal resolution of the acquired images, its spatial resolution is still limited in comparison to SMLM and STED. Better results can be achieved by Super-Resolution Radial Fluctuations (SRRF) \cite{srrf} microscopy, where a degree of local symmetry is computed for each frame. 
%While the SRRF approach has broad applicability, super-resolved SRRF images often suffer from some undesirable reconstruction artifacts. 
Recently, a SPARsity-based super-resolution COrrelation Microscopy (SPARCOM) \cite{SPARCOM} approach has been proposed to exploit the lack of correlation between different emitters as well as the sparse distribution of the fluorescent molecules via the use of a sparsity-promoting $\ell_1$-type regularization of the emitters' autocovariance. 

In the same spirit as SPARCOM, we propose a method for COvariance-based \lspace super-Resolution Microscopy with intensity Estimation (\CL) where signal sparsity is enforced in the covariance domain. Differently from SPARCOM, we use in our work an $\ell_0$-type penalization instead of the $\ell_1$ norm of the signal temporal variance. Furthermore, a joint estimation of the noise variance is performed for a more precise modeling. Our approach further incorporates a second step which allows for both signal intensity and constant background estimation, in order to deal with realistic experimental conditions. The retrieval of real image intensity information is particularly relevant for biological interpretation. 
%and for future extension to 3D super-resolution problems. 
To the best of our knowledge, \CL\ is the only super-resolution method exploiting temporal fluctuations which is capable of retrieving this information. Simulated and experimental results show that \CL\ performs well in terms of molecule localization and can retrieve accurate intensity information.

\section{\CL: formulation}
\label{sec:MethCorrCELO}

\subsection{Inverse problem formulation}
\label{ssec:model}

 Let $\mathbf{Y}_t \in \mathbb{R}^{M\times M}$ be the image frame acquired at time $t\in\left\{1,\ldots,T\right\}$, with $T>1$ and $\mathbf{X}_t \in \mathbb{R}^{L\times L}$, with $L = qM$, the desired high-resolution image defined on a $q$-times finer grid, with $q\in\mathbb{N}$.
 %and typically equal to 4. 
 The discrete model describing the acquisition process at frame $t$ is given by: 
\begin{equation} \label{eq:model}
    \mathbf{Y}_t = \mathbf{M}_q (\mathbf{H}( \mathbf{X}_t)) + \mathbf{B} +\mathbf{N_t}
\end{equation}
%the poisson distribution?    \mathbf{Y_t} = \mathbf{P}(\mathbf{M_q} (\mathbf{H}( \mathbf{X_t})) + \mathbf{B}) +\mathbf{N_t}
where $\mathbf{M}_q:\mathbb{R}^{L \times L} \rightarrow \mathbb{R}^{M \times M}$ is a down-sampling operator averaging every $q$ consecutive pixels in both dimensions, $\mathbf{H}:\mathbb{R}^{L \times L} \rightarrow \mathbb{R}^{L \times L}$ is a convolution operator defined by the PSF of the optical imaging system, $\mathbf{B}$ models the spatially and temporally constant background (due to out-of-focus fluorescent molecules) and $\mathbf{N}_t$ denotes electronic  noise, which is modeled here as a matrix of independent and identically distributed (i.i.d.) Gaussian elements, with constant variance $s \in \mathbb{R}_+$. We further assume that the molecules are located at the center of each pixel and  that there is no displacement of the specimen during the imaging period.

In vectorized form the model (\ref{eq:model}) reads:
\begin{equation} \label{eq:model_vec}
    \mathbf{y}_t = \mathbf{\Psi} \mathbf{x}_t + \mathbf{b} +\mathbf{n}_t
\end{equation}
where now $\mathbf{\Psi} \in \mathbb{R}^{M^2 \times L^2}$ is the matrix representing the composition $\mathbf{M}_q \circ \mathbf{H}$ and lower-case letters imply the column-wise vectorization of the corresponding matrices in (\ref{eq:model}). Given the matrix $\mathbf{\Psi}$ describing the acquisition process and the  frame $\mathbf{y}_t$, the task thus consists in computing a super-resolved image $\mathbf{x}$ from the many $\mathbf{y}_t$ acquisitions and, if possible, in estimating $\mathbf{b}$ and $s$. In order to exploit the statistical behaviour of the fluorescent emitters we reformulate the model in the covariance domain in the next section.

\subsection{Sparse regularization in the covariance domain}
\label{ssec:SparsityInCorDomain}

We exploit temporal and spatial independence of the fluorescent emitters by computing spatial statistics. This idea was previously exploited in \cite{sofi} and was shown to significantly reduce the full width at half maximum (FWHM) of the PSF. In particular, the use of second-order statistics 
%for a Gaussian PSF 
corresponds to a reduction factor $\sqrt2$. 
% write smth for the cross-cumulants. 

We consider the frames $\mathbf{y}_t$ as $T$ realizations of a random variable $\mathbf{y}$ whose covariance matrix is defined as:
\begin{equation}{\label{eq:cov_mat}}
    \mathbf{R_y} = \EX_{\mathbf{y}}\{(\mathbf{y} - \EX_{\mathbf{y}}\{\mathbf{y}\})(\mathbf{y} - \EX_{\mathbf{y}}\{\mathbf{y}\})^\intercal\}
\end{equation}
where $\EX_{\mathbf{y}}\{\cdot\}$ denotes the expected value computed w.r.t. to the unknown law of ${\mathbf{y}}$. We approximate $\mathbf{R_y}$ by its sample average:
\begin{equation*}{\label{eq:cov_mat_emperical}}
    \textstyle
    \mathbf{R_y} \approx
    \frac{1}{T-1}\sum_{t=1}^T (\mathbf{y}_t-\overline{\mathbf{y}})(\mathbf{y}_t-\overline{\mathbf{y}})^\intercal
    %<(\mathbf{y_t} - <\mathbf{y_t}>_t)(\mathbf{y_t} - <\mathbf{y_t}>_t)^T>_t
\end{equation*}
where $\overline{\mathbf{y}}=\frac{1}{T}\sum_{t=1}^T \mathbf{y}_t$ is the empirical temporal mean.
From (\ref{eq:model_vec}) and (\ref{eq:cov_mat}) we thus have:
\begin{equation}{\label{eq:cov_model}}
    \mathbf{R_y} = \mathbf{\Psi} \mathbf{R_x} \mathbf{\Psi}^\intercal + \mathbf{R_n} 
\end{equation}
where $\mathbf{R_x} \in \mathbb{R}^{L^2 \times L^2}$ and $\mathbf{R_n} \in \mathbb{R}^{M^2 \times M^2}$ are the covariance matrices of $\mathbf{x}_t$ and $\mathbf{n}_t$. As the background is spatio-temporally constant, the covariance matrix of $\mathbf{b}$ is zero. Recalling that the emitters are uncorrelated by assumption, we deduce that $\mathbf{R_x}$ is diagonal. We thus set $\mathbf{r_x} := \text{ diag}(\mathbf{R_x})$. Furthermore, by the i.i.d.~assumption on $\mathbf{n}_t$, we have that $\mathbf{R_n} = s \mathbf{I_{M^2}}$, where $s \in \mathbb{R}_+$ and $\mathbf{I_{M^2}}$ is the identity matrix in $\mathbb{R}^{M^2 \times M^2}$.
%Finally, due to the assumption that the background $\mathbf{b}$ is static, we have that its covariance matrix is zero. 
Equation (\ref{eq:cov_model}) can thus be re-written as:
\begin{equation*}  \label{eq:support_model}
    \mathbf{r_y} = (\mathbf{\Psi} \odot \mathbf{\Psi}) \mathbf{r_x} + s \mathbf{v_I}
\end{equation*}
where  $\odot$ denotes the column-wise Kronecker product, $\mathbf{r_y}$ is the column-wise vectorization of $\mathbf{R_y}$ and $\mathbf{v_I} = \text{vec}(\mathbf{I_{M^2}})$. 
% Setting $\Omega_{\mathbf{x}}:=\text{supp}(\mathbf{x})=\left\{i: \mathbf{x}_i\neq 0 \right\}$ and exploiting the fact that these positions are the ones of the blinking fluorescence emitters on the finer grid,  the problem thus consists in computing a sparse representation of $\Omega_{\mathbf{x}}$ and variance $s>0$ from \eqref{eq:support_model}.
In order to estimate $\mathbf{r_x} $ and the variance $s$ and promote sparsity on $\mathbf{r_x}$, we introduce the non-convex \lspace-type CEL0 regularizing penalty proposed in \cite{CELO} and consider the following minimization problem:
\begin{equation}  
    \label{eq:CEL0_variance}
    \textstyle
    \argmin\limits_{\mathbf{r_x} \geq 0, s \geq 0} \frac12 \| \mathbf{r_y} -(\mathbf\Psi \odot \mathbf\Psi) \mathbf{r_x} - s \mathbf{v_I} \|_2^2 + \Phie(\mathbf{r_x};\lambda),
\end{equation}
where $\lambda>0$ a positive regularisation parameter and the CEL0 penalty $\Phie(\cdot;\lambda)$ is defined by \cite{CELO}
\begin{equation*}
    \textstyle
    \Phie{(\mathbf{r_x};\lambda)} =  \sum\limits_{i=1}^{L^2} \lambda - \frac{\|\mathbf{a}_i\|^2}{2}\left( |(\mathbf{r_x})_i| - \frac{\sqrt{2\lambda}}{\|\mathbf{a}_i\|} \right)^2 \mathds{1} _{\{|({\mathbf{r_x}})_i| \leq \frac{\sqrt{2 \lambda}}{\|\mathbf{a}_i\|}\}},
\end{equation*}
with $\mathbf{a}_i = (\mathbf\Psi \odot \mathbf\Psi)_i$ being the i-th column of $\mathbf\Psi \odot \mathbf\Psi$. 

The functional in (\ref{eq:CEL0_variance}) is continuous and non-convex, but has the same minimizers as the corresponding $\ell_2-\ell_0$ problem. Compared to SPARCOM \cite{SPARCOM},  solutions of (\ref{eq:CEL0_variance}) are sparser. Moreover, the estimation of s is expected to improve upon the quality of the estimation. 

% Similarly to SPARCOM \cite{SPARCOM}, our approach exploits the sparsity properties of the solution $\mathbf{r_x}$  in the covariance domain, while including also the estimation of $s>0$ to improve performance. In order to promote sparsity, we make use in the following of the continuous exact relaxation of the \lspace pseudo-norm (CEL0) proposed by \cite{CELO}. Differently from the $\ell_1$-type regularisation employed in SPARCOM, the CEL0 regularisation is continuous, non-convex and preserves the minima of the $\ell_2-\ell_0$ problem one ideally would like to solve, hence its use is expected to yield sparser solutions. By further introducing some non-negativity constraints for both variables and a positive regularisation parameter $\lambda>0$, we thus aim to solve:
% \begin{equation}  \label{eq:CEL0_variance}
%     \argmin_{\mathbf{r_x} \geq 0, s \geq 0} \frac12 \| \mathbf{r_y} -(\mathbf\Psi \odot \mathbf\Psi) \mathbf{r_x} - s \mathbf{v_I} \|_2^2 + \Phie(\mathbf{r_x};\lambda),
% \end{equation}
% where $\Phie(\cdot)$ is defined by \cite{CELO}
% \begin{equation}
%     \Phie{(\mathbf{r_x};\lambda)} =  \sum\limits_{i=1}^{L^2} \lambda - \frac{\|\mathbf{a_i}\|^2}{2}\left( |(\mathbf{r_x})_i| - \frac{\sqrt{2\lambda}}{\|\mathbf{a_i}\|} \right) \mathds{1} _{\{|({\mathbf{r_x}})_i| \leq \frac{\sqrt{2 \lambda}}{\|\mathbf{a_i}\|}\}},
% \end{equation}
% with $\mathbf{a_i} = (\mathbf\Psi \odot \mathbf\Psi)_\mathbf{i}$ the i-th column of the transformation matrix. 

\subsection{Intensity and background estimation}
\label{ssec:IntensityBackground}

%As $\mathbf{r_x}$ and $\mathbf{x}$ have the same support, we deduce the support $\Omega$ of $\mathbf{x}$ from the one of $\mathbf{r_x}$. 

Solving \eqref{eq:CEL0_variance} provides an estimation of   $\mathbf{r_x}$ from which we deduce the support of $\mathbf{x}$ denoted by $\Omega_{\mathbf{x}} = \left\{i: \mathbf{x}_i\neq 0 \right\}\subset \left\{ 1,\ldots,L^2\right\}$ as $\Omega_{\mathbf{r_x}}=\Omega_{\mathbf{x}} $. We then estimate the intensity of $\mathbf{x}$ only on its support, and at the same time the spatially constant background $\mathbf{b}=b \mathbf{1_{M^2}}, b\geq 0$, by solving
%We consider a quadratic data term to model Gaussian noise 
%We are generally processing images a few non-zero elements but which are gathered and construct "continuous" structures. For this reason, we are using as a 
%and introduce a regularization term promoting smoothness on $\Omega_{\mathbf{x}}$ only, i.e.:
%Our goal is to discourage the appearance of isolated points that maybe wrongly found in the support estimation step and on the other side, encourage the 
%to retrieve smooth intensity values which are more probable to be found in real images. %For the estimation of the intensities in the subset $\Omega$ we are going to solve the following joint minimization problem:
%We thus consider the following joint minimization problem:
\begin{equation}{\label{eq:intensity}}
    \textstyle
    \argmin\limits_{\mathbf{x}\in\mathbb{R}_+^{|\Omega_{\mathbf{x}}|}, ~b \geq 0 } ~\frac12 \| \overline{\mathbf{y}} - \mathbf{\Psi_\Omega} \mathbf{x} - b \mathbf{1}_{M^2} \|_2^2 + \mu \|\nabla_{\Omega}\mathbf{x}\|_2^2
\end{equation}
%where $\mathbf{y^*} = \sum_{t=1}^T\mathbf{y_t}$ is the temporal sum of the acquired data,
where the data term models the presence of Gaussian noise, $\mu>0$ is a regularization parameter, $\mathbf{\Psi_\Omega} \in\mathbb{R}^{M^2\times |\Omega_{\mathbf{x}}|}$ is a matrix whose $i$-th column is extracted from $\mathbf\Psi$ for index $i \in \Omega_{\mathbf{x}}$ and the regularization term is the squared norm of the discrete gradient restricted to $\Omega_{\mathbf{x}}$, i.e.:
\begin{equation*}
    \textstyle
    \|\nabla_{\Omega} \mathbf{x}\|_2^2 := \frac12\sum\limits_{i \in \Omega_{\mathbf{x}}} \sum\limits_{j \in \mathcal{N}(i)\cap\Omega_{\mathbf{x}}} (x_i - x_j)^2,
\end{equation*}
where $\mathcal{N}(i)$ denotes the 8-pixel neighbourhood of $i\in\Omega_{\mathbf{x}}$. %Note that, according to this definition, $\nabla_\Omega\mathbf{x}$ denotes a (redundant) isotropic discretization of the gradient of $\mathbf{x}$ evaluated in neighbouring points for each pixel in the support $\Omega$.

%Note that this definition coincides with the one of the $\ell_2$-norm squared of a finite difference discretization of the $\nabla \mathbf{x}$ restricted to points in the support $\Omega$.

\vspace{-0.3cm}

\section{Algorithmic implementation}
\label{ssec:AlgLoc}

We use alternating minimization to solve \eqref{eq:CEL0_variance} (see Algorithm 1). For the estimation of $\mathbf{r_x}$, we follow \cite{Gazagnes} and use the iteratively reweighted $\ell_1$ algorithm (IRL1). A good initialization $\mathbf{r_x}^0$ is the $\ell_1$- regularized solution.
%which is a fast converging algorithm. % for initialization we used the solution of the l2-l1 problem  
An explicit expression for $s$ can be obtained from the  (unconstrained) optimality condition. This estimate is then projected onto the set of positive solutions in a standard way.  %which is exactly the SPARCOM solution \cite{SPARCOM}.  
For solving \eqref{eq:intensity} we use again alternating minimization (see Algorithm \ref{alg:CL2}) %as the function we aim to minimize is globally convex 
and solve each subproblem by a standard quadratic programming (QP) algorithm. Here only an initialization for the constant background value is required: a good choice for it is the median of $\overline{\mathbf{y}}$.  For both algorithms, we consider stopping criteria based on the relative difference between consecutive iterates and on a maximum number of iterations. 

%We report the pipeline in Algorithm~\ref{alg:CL1} and ~\ref{alg:CL2}. 

\vspace{-0.1cm}
\begin{algorithm}[h]
\algsetup{linenosize=\tiny}
\small{
\caption{\CL: Support estimation - 1st step}
\label{alg:CL1}
\begin{algorithmic}
\REQUIRE $\mathbf{r_y}\in\mathbb{R}^{M^4}, \mathbf{r_x}^0\in\mathbb{R}^{L^2}, s^0 \in \mathbb{R}_+, \lambda>0$
\REPEAT 
\STATE compute weights $\mathbf{\omega}_i^{\mathbf{r_x^k}}\in \partial \Phie{(\mathbf{r_x^k};\lambda)}$
\vspace{-0.5cm}
\STATE
{\begin{gather*}
\textstyle
\kern-4.2em\mathbf{r_x}^{k+1} = \argmin\limits_{\mathbf{r_x} \in \mathbb{R}_+^{L^2}} \frac12 \| \mathbf{r_y} - (\mathbf\Psi \odot \mathbf\Psi) \mathbf{r_x} - s^k \mathbf{v_I} \|_2^2
\\
\textstyle
+ \lambda \sum\limits_{i=1}^{L^2}\mathbf{\omega_i}^{\mathbf{r_x^k}}|({\mathbf{r_x}})_i|
% + \mathcal{X}_{\geq 0 }(\mathbf{r_x}) 
\end{gather*}}
\vspace{-0.2cm}
\STATE $\kern0.2em s^{k+1} = \argmin\limits_{s \in \mathbb{R}_+} \frac12 \| \mathbf{r_y} - (\mathbf\Psi \odot \mathbf\Psi) \mathbf{r_x^{k+1}} - s \mathbf{v_I} \|_2^2$
% + \mathcal{X}_{\geq 0 }(s)
\UNTIL convergence
\RETURN $\Omega_{\mathbf{x}}, s$ 
\end{algorithmic}}
\end{algorithm}
\vspace{-0.7cm}
\begin{algorithm}[h]
\algsetup{linenosize=\tiny}
\small{
\caption{\CL: Intensity estimation - 2nd step}
\label{alg:CL2}
\begin{algorithmic}
\REQUIRE $ \overline{\mathbf{y_t}}\in\mathbb{R}^{M^2}, \Omega_{\mathbf{x}}, b^0, \mu >0$
\REPEAT
\vspace{-0.6cm}
 \STATE 
{\begin{gather*}
\textstyle
\kern-1em\mathbf{x}^{k+1} = \argmin\limits_{x\in\mathbb{R}_+^{|\Omega_{\mathbf{x}}|}} \frac12 \|  \overline{\mathbf{y}} - \mathbf{\Psi_\Omega} \mathbf{x} - b^{k} \mathbf{1_{M^2}} \|_2^2 + \mu \|\nabla_{\Omega}\mathbf{x}\|_2^2\\
\vspace{-0.3cm}
%\emph{\quad{\% By QP}}
%+ \mathcal{X}_{\geq 0 }(\mathbf{x})
\end{gather*}}
\vspace{-1cm}
\STATE $\kern0.2em b^{k+1} = \argmin\limits_{b \in \mathbb{R}_+} \frac12 \| \overline{\mathbf{y}} - \mathbf{\Psi_\Omega} \mathbf{x}^{k+1} - b \mathbf{1_{M^2}} \|_2^2 $
%+ \mathcal{X}_{\geq 0 }(b)
\UNTIL convergence
\RETURN $\mathbf{x}, b$
\end{algorithmic}}
\end{algorithm}

\vspace{-0.4cm}

\section{Results and discussion}
\label{sec:Results}

\subsection{Simulated data}
\label{ssec:gensimulated}
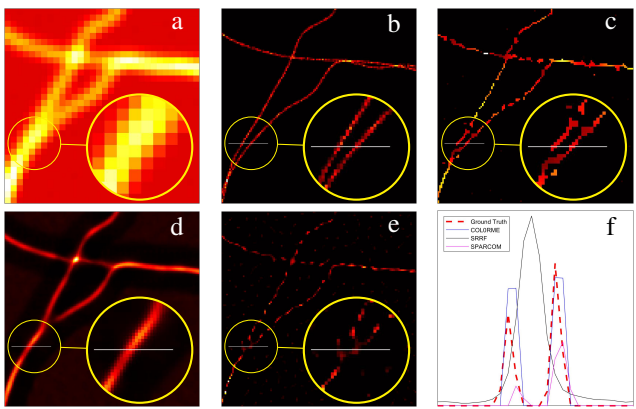
\begin{figure}[!htbp]%
    \captionsetup[subfigure]{labelformat=empty}
    \vspace{-0.5cm}
    \centering
    \subfloat[\centering]{\begin{tikzpicture}[spy using outlines={circle,yellow,magnification=2,size=1.4cm, connect spies}]
        \node {\pgfimage[width=2.6cm]{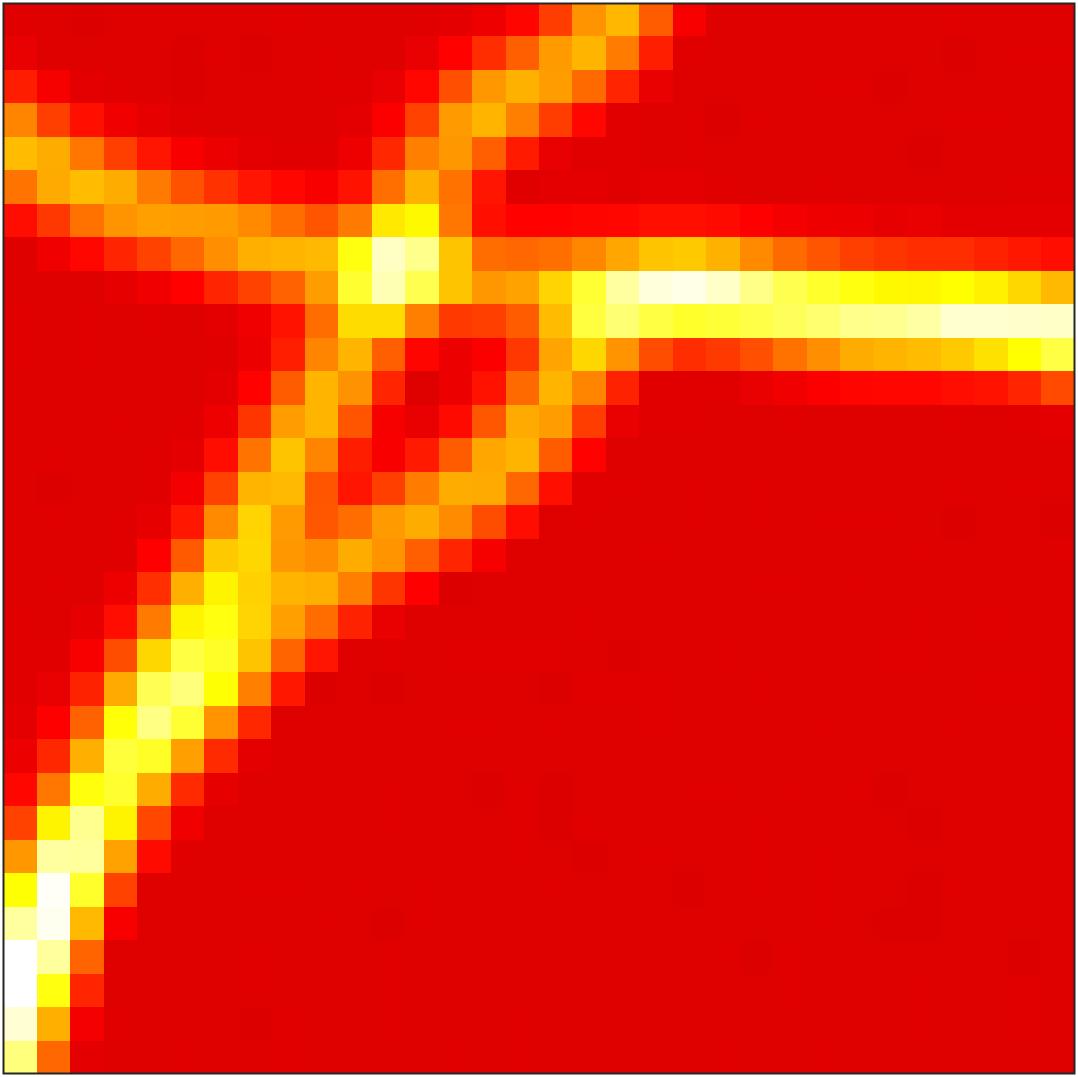}};
        \spy on (-0.9,-0.5) in node [left] at (1.2,-0.55);
        \node at (1,1.1) {\color{white}a};
        \end{tikzpicture}}
    \hfill
    \subfloat[\centering]{\begin{tikzpicture}[spy using outlines={circle,yellow,magnification=2,size=1.4cm, connect spies}]
        \node {\pgfimage[width=2.6cm]{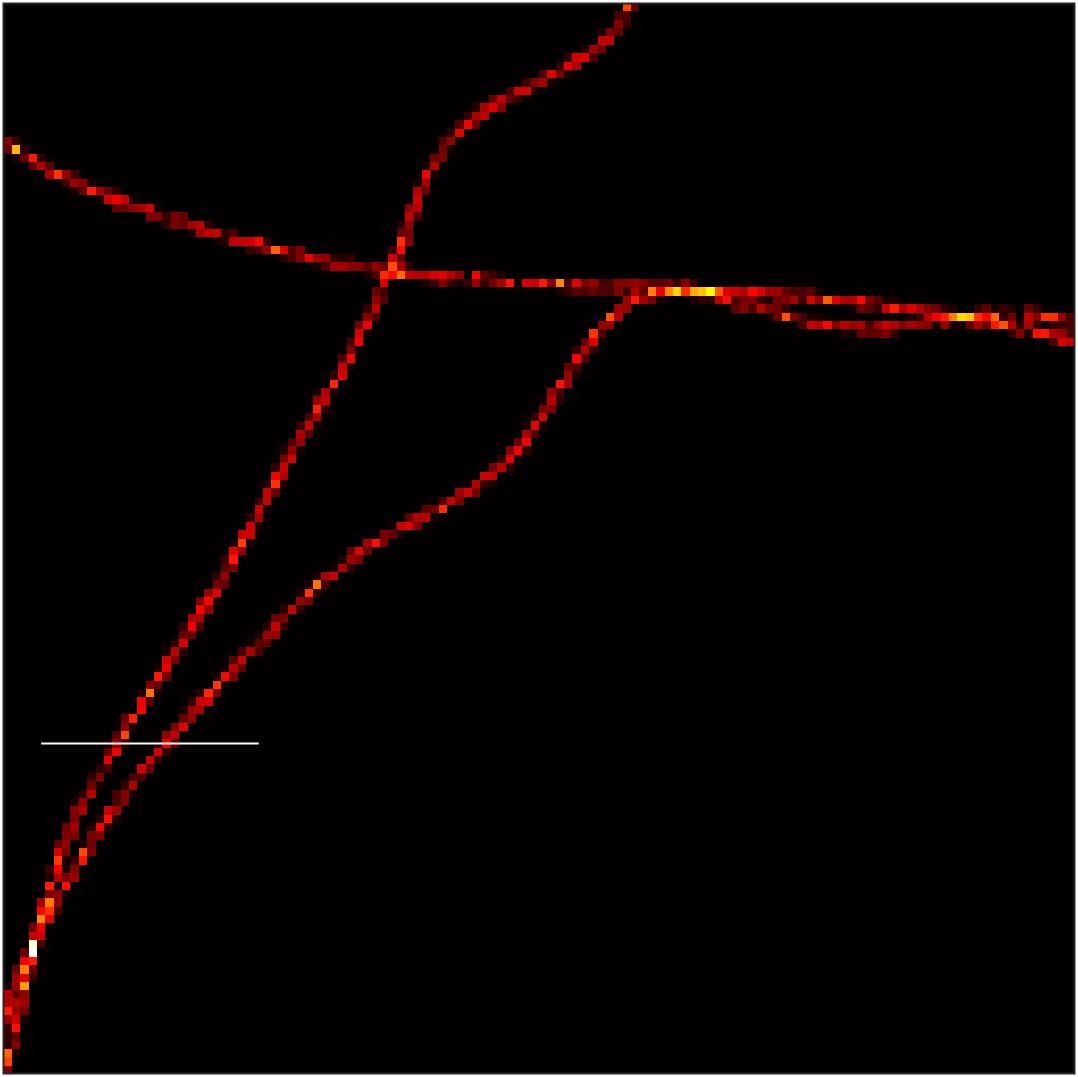}};
        \spy on (-0.9,-0.5) in node [left] at (1.2,-0.55);
        \node at (1,1.1) {\color{white}b};
        \end{tikzpicture}\label{fig:gt}}
    \hfill
    \subfloat[\centering]{\begin{tikzpicture}[spy using outlines={circle,yellow,magnification=2,size=1.4cm, connect spies}]
        \node {\pgfimage[width=2.6cm]{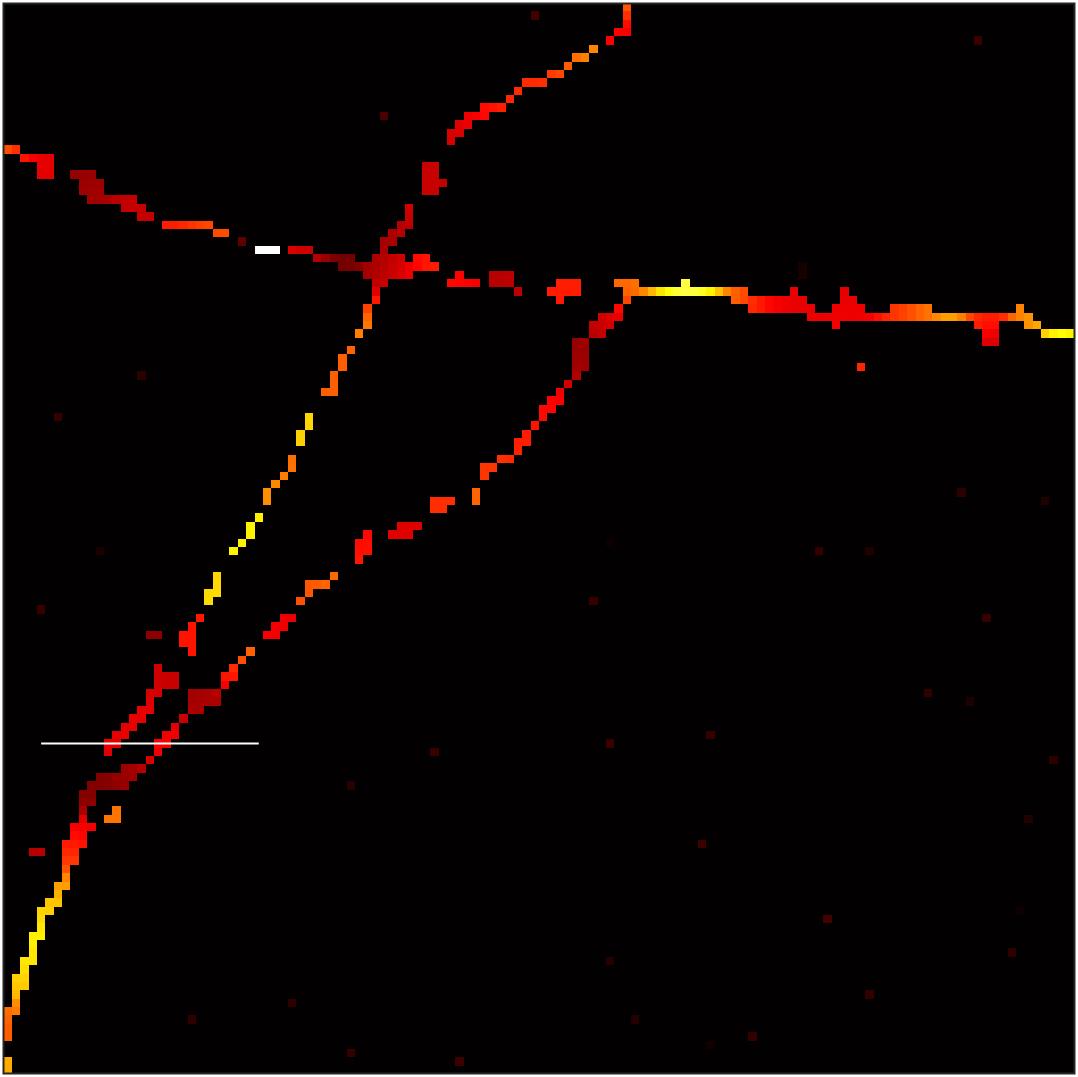}};
        \spy on (-0.9,-0.5) in node [left] at (1.2,-0.55);
        \node at (1,1.1) {\color{white}c};
        \end{tikzpicture}}
    \vspace{-1cm}
    \centering
    \subfloat[\centering]{\begin{tikzpicture}[spy using outlines={circle,yellow,magnification=2,size=1.4cm, connect spies}]
        \node {\pgfimage[width=2.6cm]{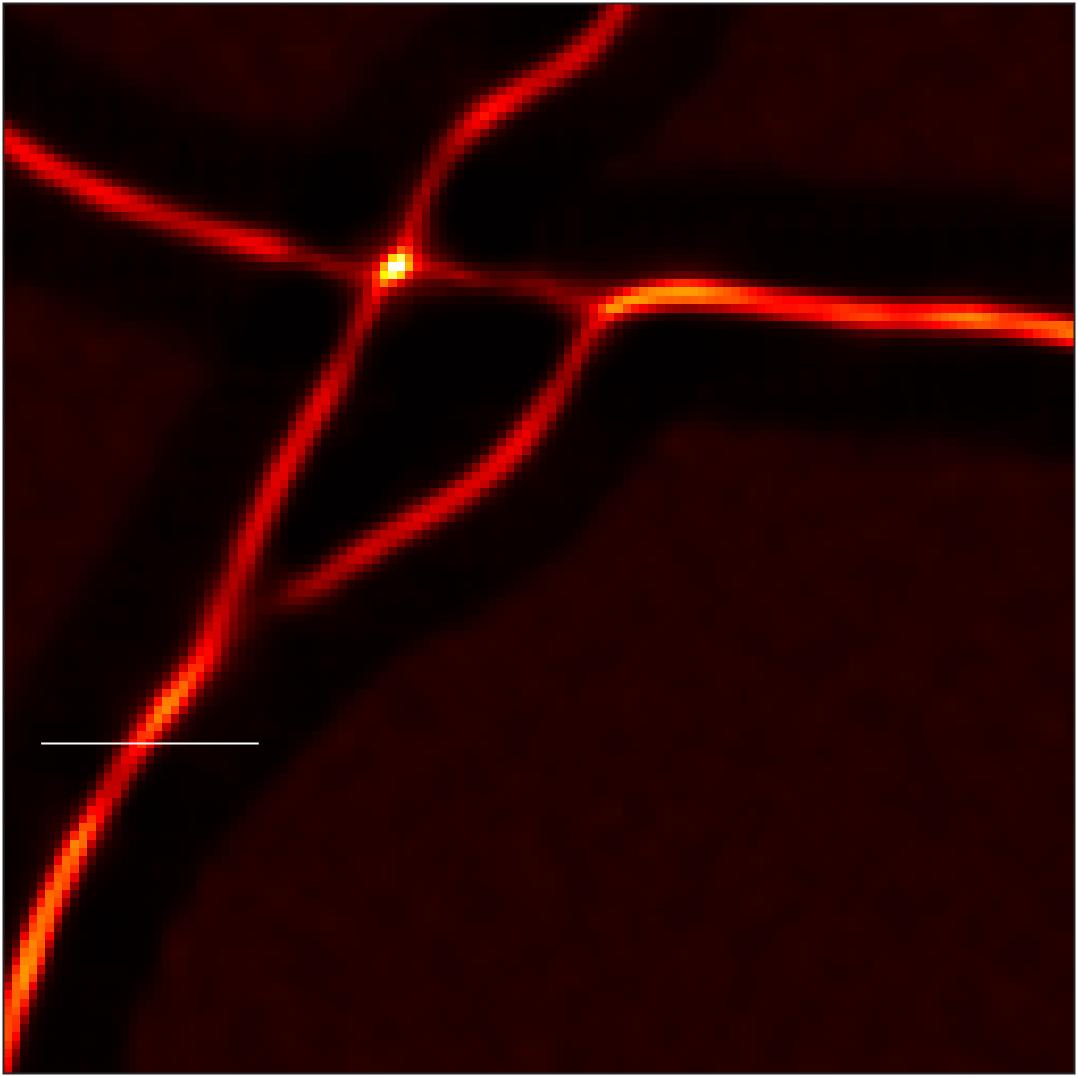}};
        \spy on (-0.9,-0.5) in node [left] at (1.2,-0.55);
        \node at (1,1.1) {\color{white}d};
        \end{tikzpicture}}
    \hfill
    \subfloat[\centering]{\begin{tikzpicture}[spy using outlines={circle,yellow,magnification=2,size=1.4cm, connect spies}]
        \node {\pgfimage[width=2.6cm]{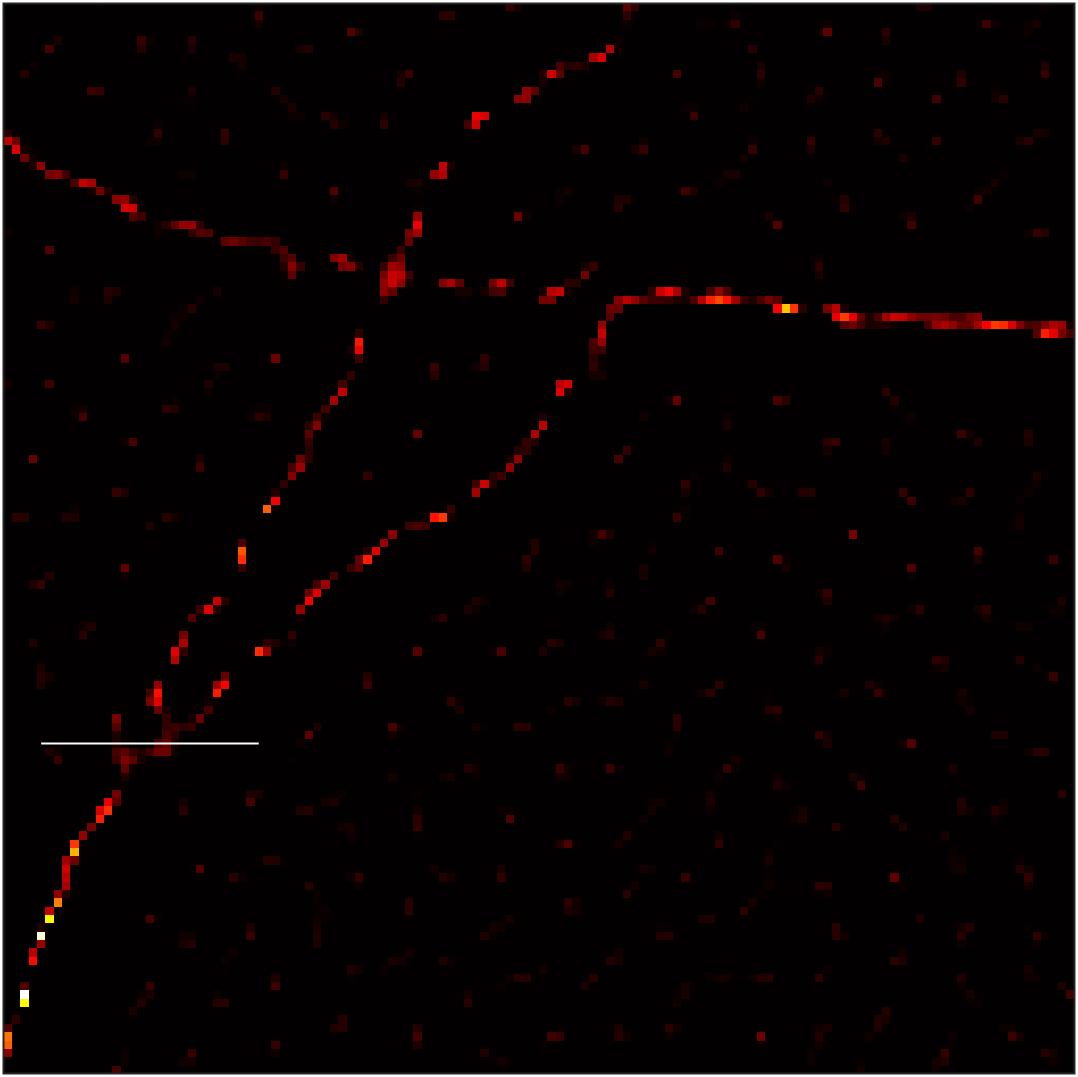}};
        \spy on (-0.9,-0.5) in node [left] at (1.2,-0.55);
        \node at (1,1.1) {\color{white}e};
        \end{tikzpicture}}
    \hfill
    \subfloat[\centering]{\begin{tikzpicture}[spy using outlines={circle,yellow,magnification=2,size=1.4cm, connect spies}]
        \node {\pgfimage[width=2.6cm]{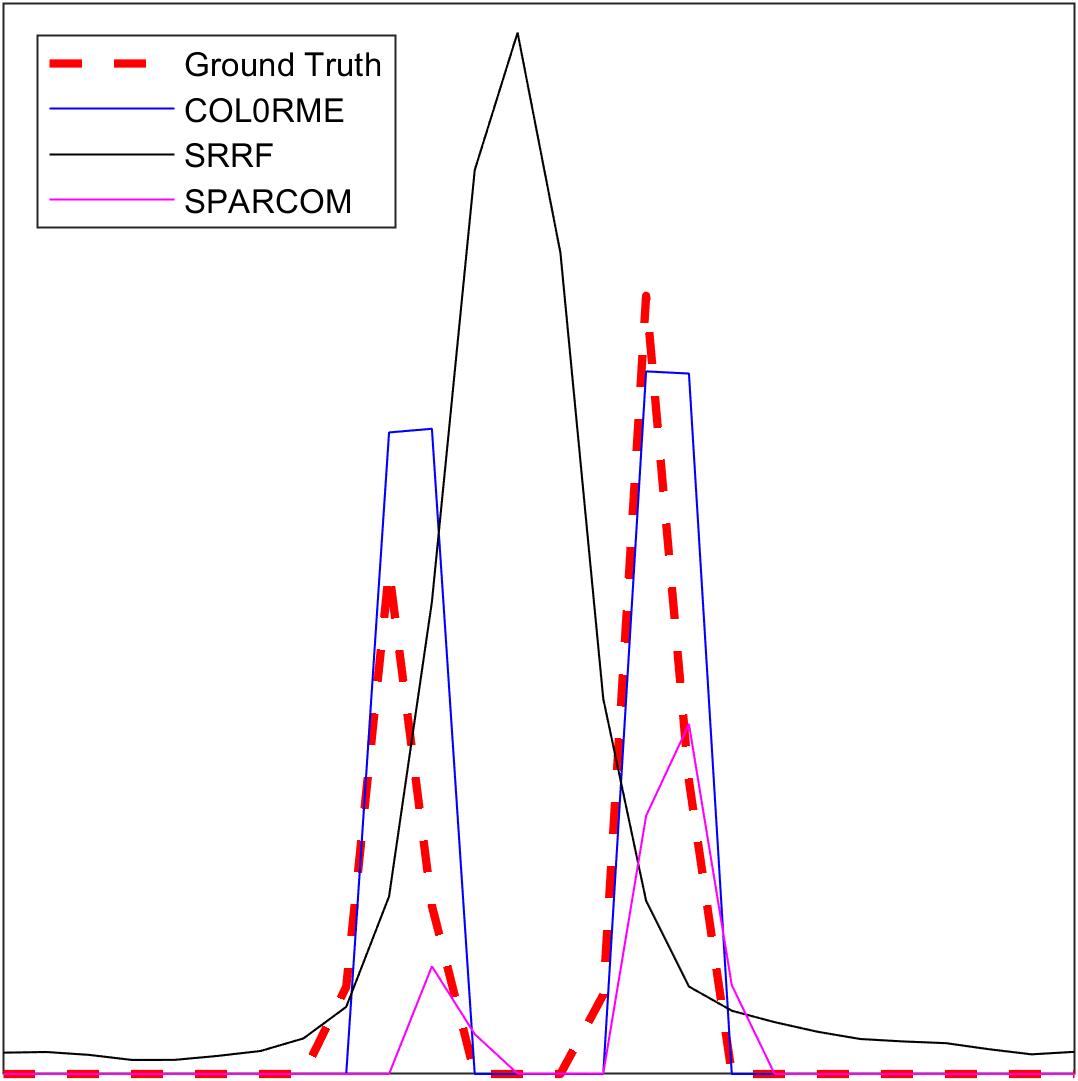}};
        \node at (1,1.1) {f};
        \end{tikzpicture}\label{fig:lineplot}}
    \vspace{-0.7cm}
    \caption{ \small{Results for the simulated 'High Background' dataset,  with $T=700$ before PP: (a) temporal average of the acquired stack ($4\times$ zoom), (b) ground truth, (c) \CL, (d) SRRF, (e) SPARCOM, (f) intensity profiles (SRRF and SPARCOM with adapted range)}}%
    \label{fig:reconstructions}
    \vspace{-0.5cm}
\end{figure}

We start by applying \CL\ to images of tubular structures simulating standard microscope acquisitions with standard fluorescent dyes. The spatial pattern is taken from the MT0 microtubules training dataset uploaded for the SMLM  2016\footnote{\href{http://bigwww.epfl.ch/smlm/datasets/index.html}{http://bigwww.epfl.ch/smlm/datasets/index.html}}, see Fig. \ref{fig:gt}. Intensities are obtained by using the SOFI simulation tool \cite{SOFItool}. Namely, we simulate temporal fluctuations and create videos of $T=100$ and $T=700$ frames at a frame rate of 100 frames per second (fps). The simulation parameters are set as follows: $20$ms for on-state average lifetime, $40$ms for off-state average lifetime and $20$s for average time untill bleaching.  The emitter density is equal to 10.7 emitters/pixel/frame, while the FWHM of the PSF is approximately $229$nm.  We create two different noisy datasets in order to evaluate the results of \CL\ and to compare them with the ones obtained by SRRF\cite{srrf} and SPARCOM\cite{SPARCOM}, which, similarly, also exploit the temporal fluctuations of molecules. For the first dataset, we generate on average 1000 photons/frame per emitting molecule and $b^*=100$ photons/frame per pixel to simulate the out-of-focus molecules, which we consider as background (BG). For the second dataset (which is more realistic), we set these values to 500 and $b^*=2500$, respectively. For both datasets we add Gaussian noise of $20$dB to simulate the presence of electronic noise.
%The signal-to-noise ratio (SNR) in the first dataset is  21.26 dB, while in the second (which is more realistic) it is 19.64 dB. 

As we know the ground truth positions of the emitters, we can evaluate the localization precision of each method in terms of the Jaccard index, a quantity in the range $[0,1]$ computed as the ratio between correct detections (up to some tolerance $\delta>0$) and the sum of correct detections, false positives and false negatives. Table \ref{table:Jacc} reports the average Jaccard index ($\delta=40$nm) for 20 noise realizations. Note that, in this table, for both SRRF and SPARCOM the reported results have been obtained after a thresholding step to eliminate the many false positives in the reconstructed support, while \CL\, does not require any post-processing and achieves precise localization even in hard experimental conditions. In Fig. \ref{fig:reconstructions} we report the final reconstructions  before thresholding showing that \CL\, better preserves the real structure of the tubulin. This is a different feature from SRRF whose reconstruction is overall visually more pleasant, but which creates artifacts as it cannot separate close-by filaments. \CL\ exhibits  a more accurate reconstruction than SPARCOM whose result has many false positives around the tubulin. Thanks to the intensity estimation step, \CL\, is also able to estimate the brightness values which correspond to real signal intensities. The peak signal-to-noise ratio (PSNR) values w.r.t. to the ground truth intensity image and the estimated BG value $b$ are provided in Table \ref{table:PSNR}. We remark that the SRRF and SPARCOM do not actually estimate real signal intensities. Consequently, we had to adapt the range of their reconstructions according to the range of the ground truth image for the line profile graph in Fig. \ref{fig:lineplot}. 

\begin{table}[!htb]
\centering
\begin{tabular}{l|c|c|c|c|}
\cline{2-5}
 & \multicolumn{2}{c|}{Low BG} & \multicolumn{2}{c|}{High BG} \\ \hline
\multicolumn{1}{|c|}{\backslashbox{Method}{$T$}} & 100 & 700 & 100 & 700 \\ \Xhline{3\arrayrulewidth}
\multicolumn{1}{|l|}{\CL~ (No PP)}                       & \textbf{0.51} & \textbf{0.66} & 0.28 & \textbf{0.51} \\ \hline
\multicolumn{1}{|l|}{SRRF (PP)}                 & 0.5  & 0.5 & \textbf{0.47} & 0.47 \\ \hline
\multicolumn{1}{|l|}{SPARCOM (PP)}              & 0.32  & 0.39  & 0.17 & 0.31 \\ \hline
\end{tabular}%
\vspace{-0.1cm}
\caption{\small{Jaccard index values with tolerance $\delta = 40$nm for the different datasets. The results obtained by SRRF and SPARCOM are subject to a post-processing (PP) step to remove false positives. 
%No PP is needed for \CL.}
}}
\label{table:Jacc}
\end{table}

\begin{table}[!htb]
\begin{tabular}{lllll}
\cline{2-5}
\multicolumn{1}{l|}{}      & \multicolumn{2}{c|}{Low BG, $b^*=100$}           & \multicolumn{2}{c|}{High BG, $b^*=2500$}          \\ \hline
\multicolumn{1}{|l|}{$T$}    & \multicolumn{1}{c|}{100} & \multicolumn{1}{c|}{700} & \multicolumn{1}{c|}{100} & \multicolumn{1}{c|}{700} \\ \Xhline{3\arrayrulewidth}
\multicolumn{1}{|l|}{PSNR} & \multicolumn{1}{c|}{$26.6$} & \multicolumn{1}{c|}{$29.21$} & \multicolumn{1}{c|}{$23.31$} & \multicolumn{1}{c|}{ $28.37$} \\ \hline
\multicolumn{1}{|l|}{$b$} & \multicolumn{1}{c|}{$92.3$} & \multicolumn{1}{c|}{$102$} & \multicolumn{1}{c|}{$2435$} & \multicolumn{1}{c|}{ $2454$} \\
\hline
\end{tabular}
\caption{\small{\CL\ PSNR values (in dB) and estimated background $b$ for the four datasets (average over 20 noise realizations).}}
\label{table:PSNR}
\end{table}

\vspace{-0.6cm}

\vspace{0.5cm}
\subsection{Experimental data}
\label{ssec:experimental}
\begin{figure}[!tbp]%
    \centering
    \captionsetup[subfigure]{labelformat=empty}
    \subfloat[\centering]{\begin{tikzpicture}
        \node {\pgfimage[width=2.6cm]{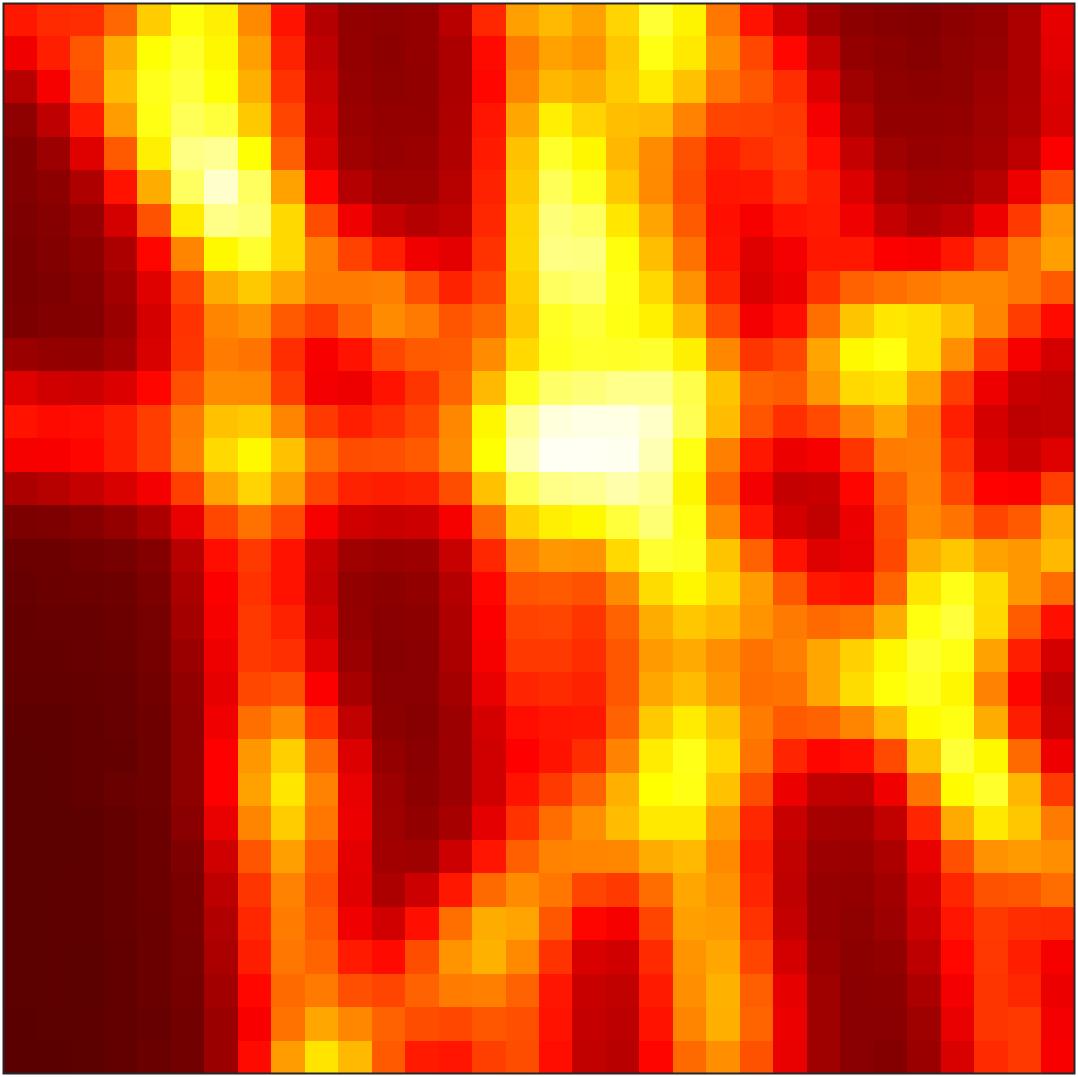}};
        \node at (1,1.1) {\color{white}a};
        \end{tikzpicture}}
    \hfill
    \subfloat[\centering]{\begin{tikzpicture}
        \node {\pgfimage[width=2.6cm]{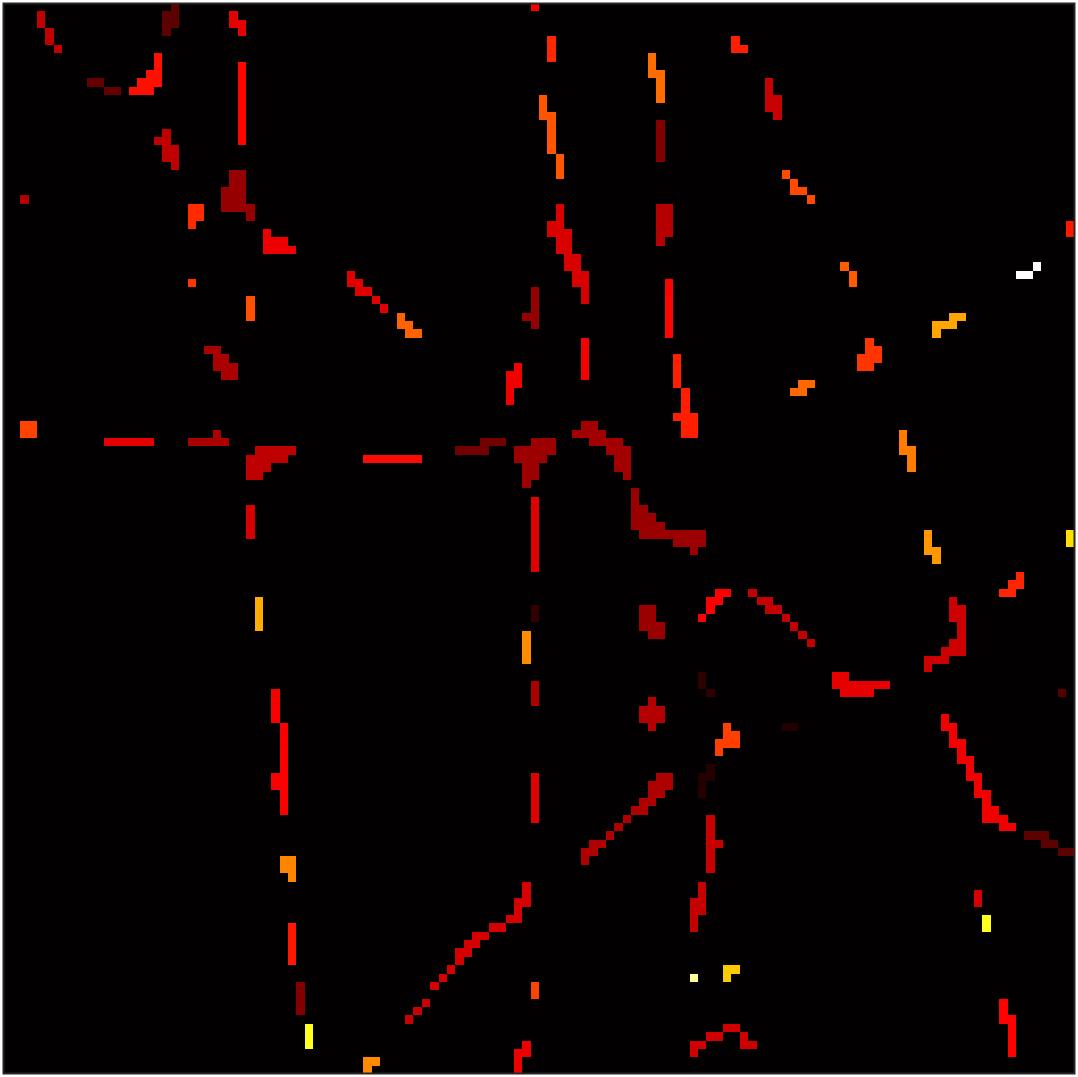}};
        \node at (1,1.1) {\color{white}b};
        \end{tikzpicture}}
    \hfill
    \subfloat[\centering]{\begin{tikzpicture}
        \node {\pgfimage[width=2.6cm]{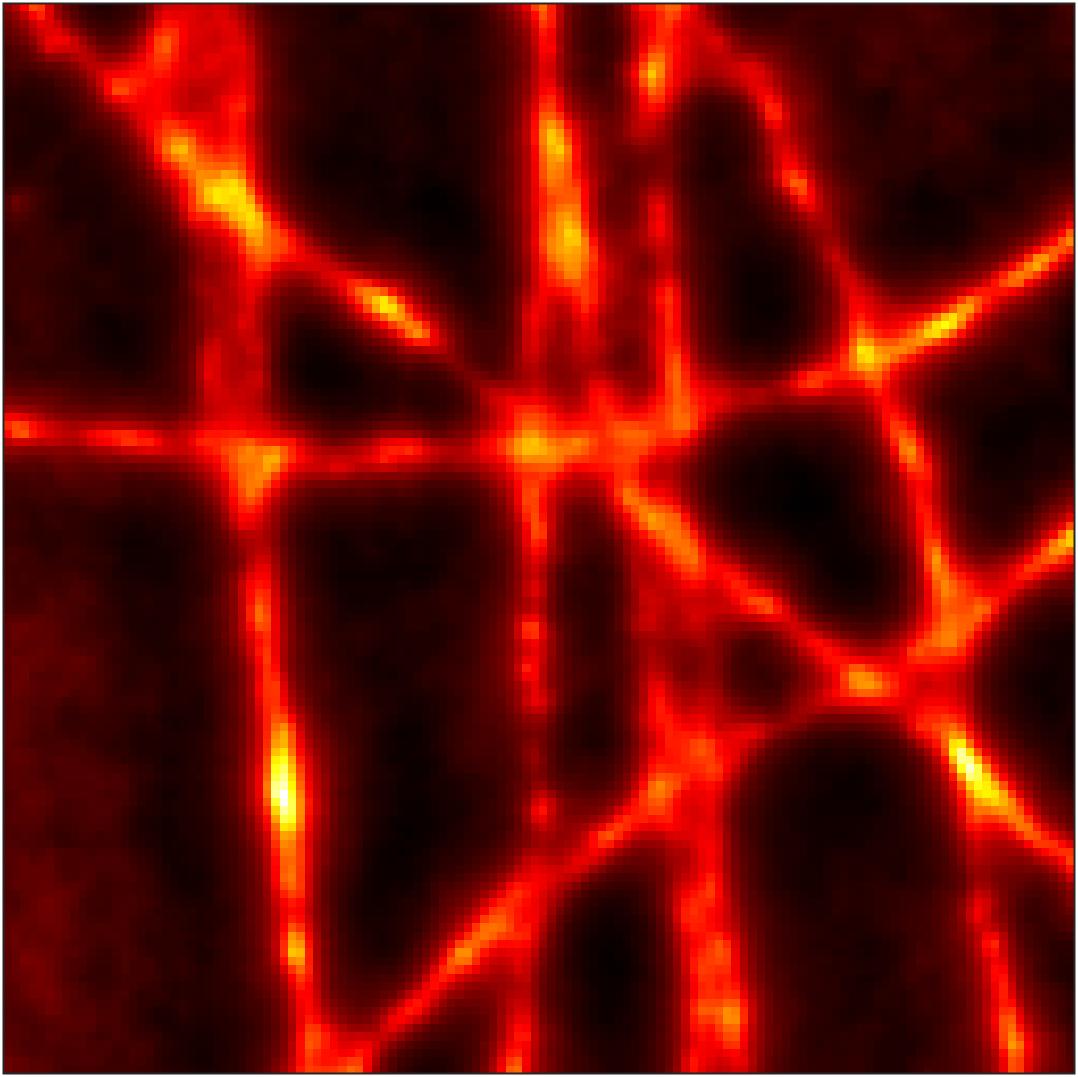}};
        \node at (1,1.1) {\color{white}c};
        \end{tikzpicture}}
    \vspace{-0.7cm}
    \caption{\small{ Results on high-density SMLM data: (a) summation of the acquired stack ($4\times$ zoom), (b) \CL, (c) SRRF}}%
    \label{fig:reconstructionsExp}
    \vspace{-0.4cm}
\end{figure}
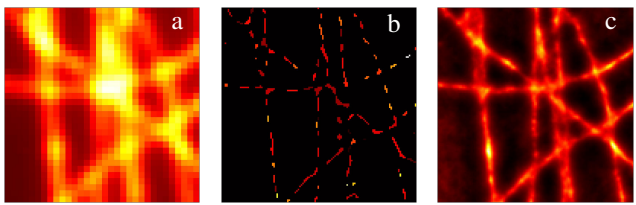

\vspace{-0.05cm}
\CL~ can also be applied  for high-density acquisitions obtained by SMLM techniques. Even though in SMLM the molecules do not have a blinking behaviour, but rather a on-to-off transition, we can consider as blinking the temporal behaviour of one pixel in high density videos with many molecules per pixel. We compared the methods \CL\, and SRRF on a patch extracted from a real dataset from the SMLM challenge 2013. The acquired video has $T=500$ frames and a FWHM of $351.8$ nm. In Fig. \ref{fig:reconstructionsExp} we observe that SRRF preserves better the broad structure of the specimen, while \CL\ better separates very close microtubules and does not create any background artifacts. %The reconstructed structures are not continuous in this example but the localization is more precise and there are no artifacts as in the method SRRF. 

\section{Conclusions}
\label{sec:conclusions}
\vspace{-0.20cm}
We proposed a super-resolution method named \CL, well suited for live-cell imaging and which can be easily applied to images obtained by common microscopes and conventional dyes.\;\CL\; takes advantage of the independent stochastic fluctuations of the fluorescent molecules and by solving a non-convex optimization problem in the covariance domain localizes the fluorophores with high precision. Differently from other super-resolution methods (SRRF, SPARCOM), \CL~ also includes an intensity estimation step, which is a valuable piece of information for biological analysis. We showed that \CL\; outperforms competing methods in terms of localization precision on both simulated and real data, and further computes an estimate of noise statistics and background information.   
\vspace{-0.13cm}

\small{\section{Compliance with ethical standards}

\vspace{-0.25cm}

This work was conducted using biological data available in open access by EPFL SMLM datasets. Ethical approval was not required as confirmed by the license attached with the open access data.}
\vspace{-0.13cm}

\small{\section{Acknowledgments}

\vspace{-0.25cm}

The work of VS and LBF has been supported by the French government, through the 3IA Côte d’Azur Investments in the Future project managed by the National Research Agency (ANR) with the reference number ANR-19-P3IA-0002. LC acknowledges the support received by
UCA IDEX JEDI and by the NoMADS RISE H2020 project 777826.}

\normalsize
\vspace{-0.25cm}

\bibliographystyle{IEEEbib}
\bibliography{strings}

\end{document}